\documentclass[12pt]{amsart}

\usepackage[usenames]{color}

\usepackage{amsmath}
\usepackage{amsthm}
\usepackage{amsfonts}
\usepackage{amssymb}
\usepackage{latexsym}
\usepackage{color}

\usepackage{hyperref}
\usepackage[top=1in, bottom=1in, left=1in, right=1in]{geometry}

\newtheorem{theorem}{Theorem}
\newtheorem{corollary}{Corollary}

\newtheorem{proposition}{Proposition}

\newcommand{\Z}{\mathbb{Z}}

\newcommand{\F}{\mathbb{F}}

\begin{document}

\author[]{Nadya Nabahi}
\address{Department of Mathematics, The City College of New York, New York,
NY 10031} \email{nnabahi000@citymail.cuny.edu}

\author[]{Vladimir Shpilrain}
\address{Department of Mathematics, The City College of New York, New York,
NY 10031} \email{shpilrain@yahoo.com}

\title{Easy estimates of Lyapunov exponents\\
for random products of matrices}

\begin{abstract}
The problems that we consider in this paper are as follows. Let $A_1, \ldots, A_k$ be square  matrices (over reals). Let $W=w(A_1, \ldots, A_k)$ be a random product of $n$ matrices. What is the expected growth rate of the largest (in the absolute value) entry in such a random product? What is the (maximal) Lyapunov exponent for a random matrix product like that? We give  an answer to the first question under some mild restrictions on the entries of $A_i$. For the second question, we offer a very simple and efficient method to produce an upper bound on the Lyapunov exponent.

\end{abstract}

\maketitle

\noindent {\it Keywords:} random matrix product, Lyapunov exponent

\noindent {\it MSC numbers:} 15B52, 60B20

\section{Introduction}

The first problem that we consider in this paper is as follows. Let $A_1, \ldots, A_k$ be square  matrices (over reals). Let $W=w(A_1, \ldots, A_k)$ be a random product of $n$ matrices. What is the expected growth rate of the largest (in the absolute value) entry in such a random product? 
Or, equivalently: what is the expected growth rate of the norm of $W$?

To make the exposition easier to follow, we do what many other authors do: we limit the exposition to the case of just two matrices, $A$ and $B$. By a ``random product" of $A$ and $B$ we will typically mean a product of matrices where each matrix is either $A$ (with probability $p$) or $B$ (with probability $1-p$).

We show that the above problem has a surprisingly simple solution:

\begin{proposition} \label{thm1}
\noindent {\bf (a)} Suppose all entries of the matrices $A$ and $B$ are nonnegative real numbers. Suppose also that all eigenvalues of the matrix $M= pA + (1-p)B$ are real numbers. Let the largest in the absolute value eigenvalue of $M$, call it $\mu$, be positive and unique (i.e., there is no other eigenvalue of $M$ with the same absolute value). Then the expectation of the largest entry of a random product $w(A, B)$ of length $n$ is $\Theta(\mu^n)$. That is, the expected growth rate of the largest (in the absolute value) entry in $w(A, B)$ is precisely $\mu$.
\medskip

\noindent {\bf (b)} For a matrix $S$, denote by $|S|$ the matrix where all entries of $S$ are replaced by their absolute values. Let $\bar{M}=p|A| + (1-p)|B|$.  Let the largest in the absolute value eigenvalue of $\bar{M}$, call it $\bar{\mu}$, be positive and unique.
Then the expectation of the largest (in the absolute value) entry of a random product $w(A, B)$ of length $n$ is $O(\bar{\mu}^n)$.
\end{proposition}

Proposition \ref{thm1} applies, in particular, to the case where matrices are $2 \times 2$, and the matrix
$M= pA + (1-p)B$ has only positive entries; then the existence of the largest (positive) real eigenvalue (the Perron–Frobenius eigenvalue) is guaranteed by the Perron–Frobenius theorem. However, conditions of Proposition \ref{thm1} are satisfied for many matrices with negative entries as well.

The proof of Proposition \ref{thm1} is given in Section \ref{growthav}.

A more difficult problem, of interest in the theory of stochastic processes, is to estimate
$L=\lim_{n\to \infty} \sqrt[n]{||w(A, B)||}$, where $||w(A, B)||$ denotes the norm of the matrix $W=w(A, B)$. The natural logarithm of $L$ is what is called the (maximal) Lyapunov exponent of a random product of $A$ and $B$, see e.g. \cite{Pollicott}. We note that one of the ways to define the norm of a matrix is to take its largest (in the absolute value) entry, see Section \ref{proof2}.

For strictly positive non-singular matrices, Pollicott \cite{Pollicott} reported an algorithm to estimate the Lyapunov exponent  $\lambda = \log L$ with any desired precision.
Pollicott mentions that ``it is a fundamental problem to find both an explicit expression for $\lambda$ and a useful method of accurate approximation". See also \cite{Tsitsiklis}.

Our contribution to this problem is a very simple and efficient method to produce an upper bound on the Lyapunov exponent based on Proposition \ref{thm1} and the following

\begin{theorem} \label{thm2}
Let $A$ and $B$ be $d \times d$ matrices over reals.
Suppose the matrices $A$, $B$, and $\bar{M}$, where $\bar{M}=p|A| + (1-p)|B|$, satisfy conditions of Proposition \ref{thm1}(b),
and let $\mu$ be the (unique) largest positive eigenvalue of $M'$. Then the Lyapunov exponent $\lambda$ of a random product of matrices $A$ and $B$ satisfies $\lambda \le \log \mu$ with probability 1.

\end{theorem}

Proposition \ref{thm1} and Theorem \ref{thm2} generalize to an arbitrary finite collection $\{A_1, \ldots, A_k\}$ of matrices, not just for pairs. 

Theorem \ref{thm2}, in combination with Proposition \ref{thm1}, gives a very simple and efficient method for bounding the Lyapunov exponent from above. We note that altogether different (analytical) methods for bounding the Lyapunov exponent were previously used in the literature, see e.g. \cite{Jurga}, \cite{Pollicott2}, \cite{Protasov}, and \cite{Sturman}. We compare our bounds to theirs for some particular pairs of matrices in Sections \ref{upper} and \ref{Examples}, and summarize advantages and disadvantages of our method in the concluding Section \ref{Conclusions}.

We also mention that our method applies to a larger class of matrices than methods of \cite{Jurga}, \cite{Pollicott2}, \cite{Sturman} and others do. In particular, we do not require all the entries of participating matrices $A_i$ to be strictly positive. Neither do we have any kinds of non-singularity conditions on matrices $A_i$. For example, in \cite{Pollicott} one of the requirements is that all matrices $A_i$ are non-singular. In \cite{Protasov}, this was relaxed to $A_i$ not having a row or a column of zeros. With our approach, if we take, say, $A = \left(
 \begin{array}{cc} 0 & 0 \\ 2 & 1 \end{array} \right), B = \left(
 \begin{array}{cc} 0 & 2 \\ 0 & 1 \end{array} \right)$, $p=\frac{1}{2},$ then  $M=pA + (1-p)B = \left( \begin{array}{cc} 0 & 1 \\ 1 & 2 \end{array} \right)$, and this matrix $M$ satisfies all conditions of Theorem \ref{thm2}.
\medskip

Finally, we note that it is also an important problem to determine the {\it maximum possible} (in the absolute value) entry $G(n)$ in a product $w(A_1, \ldots, A_k)$ of length $n$. The number $\lim_{n \to \infty}\sqrt[n]{|G(n)|}$ (i.e., the maximum growth rate) is called the {\it joint spectral radius} $\rho(A_1, \ldots, A_k)$ of the collection $\{A_1, \ldots, A_k\}$ of matrices, see e.g. \cite{Jungers}. Computing the joint spectral radius is a difficult problem in general, although there are results in various special cases, e.g. \cite{BSV}, \cite{Masuda}, \cite{Lagarias}, \cite{Panti}. Computing the joint spectral radius of pairs of matrices over integers has important real-life applications, see \cite{survey} for a survey.

Since the average growth rate is obviously not larger than the maximum growth rate, Proposition \ref{thm1}(a) implies:

\begin{corollary}\label{spectral}
Suppose the matrices $A$,  $B$, and $M=pA + (1-p)B$ satisfy conditions of Proposition \ref{thm1}(a),
and let $\mu$ be the (unique) largest positive eigenvalue of $M$. Then the joint spectral radius
$\rho(A, B)$ satisfies the inequality $\rho(A, B) \ge \mu$.

\end{corollary}

Just like Proposition \ref{thm1}(a), Corollary \ref{spectral} generalizes to arbitrary collections  $\{A_1, \ldots, A_k\}$ of matrices, not just pairs.

\section{Average Growth}\label{growthav}

Perhaps surprisingly, computing the {\it average} growth rate of the entries in a random product of $n$ matrices $A$ and $B$ (where each factor is $A$ or $B$ with probability $p$ and $1-p$, respectively) is easier than computing the {\it maximum} growth rate (i.e., the joint spectral radius).

A straightforward method for computing the average growth rate is based on solving a system of linear recurrence relations with constant coefficients. To illustrate this method, we start with
a couple of simple examples. In what follows, we denote $A(k) = \left(
 \begin{array}{cc} 1 & k \\ 0 & 1 \end{array} \right), \hskip .2cm B(m) = \left(
 \begin{array}{cc} 1 & 0 \\ m & 1 \end{array} \right)$.

\medskip

\subsection{Average growth for products of $A(1)$ and $B(1)$}\label{growthav11}\hfill \break

\noindent  Let $\left(\begin{array}{cc} a_n & b_n \\ c_n & d_n \end{array} \right)$ denote the result of multiplying $n$ matrices where each factor is $A=A(1)$ or $B=B(1)$ with probability $\frac{1}{2}$. Denote the expectation of $a_n$ by  $\bar a_n$, etc. Then, using linearity of the expectation, we have the following recurrence relations for the expectations. (Note that since all entries of the matrices involved are nonnegative, we have in this case $|\bar a_n| = \bar a_n$ and $|\bar b_n| = \bar b_n$ for all $n$.)
\medskip

\noindent $|\bar a_n| = (|\bar a_{n-1}|, ~|\bar b_{n-1}|)\cdot \left(\frac{1}{2} \left(\begin{array}{c} 1 \\0\end{array} \right) + \frac{1}{2} \left(\begin{array}{c} 1 \\1\end{array} \right)\right) =
\frac{1}{2} |\bar a_{n-1}| + \frac{1}{2}(|\bar a_{n-1}| + |\bar b_{n-1}|) = |\bar a_{n-1}| + \frac{1}{2} |\bar b_{n-1}|$.
\medskip

\noindent $ |\bar b_n| = (|\bar a_{n-1}|, ~|\bar b_{n-1}|)\cdot \left(\frac{1}{2} \left(\begin{array}{c} 1 \\1\end{array} \right) + \frac{1}{2} \left(\begin{array}{c} 0 \\1\end{array} \right)\right) =
\frac{1}{2}(|\bar a_{n-1}| + |\bar b_{n-1}|) + \frac{1}{2} |\bar b_{n-1}| = \frac{1}{2} |\bar a_{n-1}| + |\bar b_{n-1}|$.
\medskip

\noindent The same recurrence relations hold (independently) for $|\bar c_n|$ and $|\bar d_n|$, so it is sufficient to handle just $|\bar a_n|$ and $|\bar b_n|$.

Following the usual method of solving a ``square" system of linear recurrence relations with constant coefficients, we write it in the matrix form:
\medskip

\noindent $(|\bar a_n|, |\bar b_n|) = M \cdot (|\bar a_{n-1}|, |\bar b_{n-1}|)$, where vectors on the left and on the right should be interpreted as columns. The matrix $M$ here is $M = \left( \begin{array}{cc} 1 & \frac{1}{2} \\ \frac{1}{2} & 1 \end{array} \right) = \frac{1}{2} A + \frac{1}{2} B$.

The solution of the system therefore is $(|\bar a_n|, |\bar b_n|) = M^n \cdot (1, 0)$.
The largest eigenvalue of the matrix $M$ in this case is $\frac{3}{2}$, so the average growth rate is $\mu=\frac{3}{2}$.
\medskip

\subsection{Average growth for products of $A$ and $B$ if negative entries are allowed}\label{growthav1-1}\hfill \break

Now let us see what happens if negative entries in matrices $A$ and/or $B$ are allowed, e.g. if $A=\left(\begin{array}{cc} -2 & 1 \\ 3 & 1 \end{array} \right)$ and $B= \left(\begin{array}{cc} 1 & 0 \\ 1 & 1 \end{array} \right)$.
The recurrence relations then are:

\medskip

\noindent $|\bar a_n| = \frac{1}{2}  |(|\bar a_{n-1}|, |\bar b_{n-1}|) \cdot \left(\begin{array}{c} -2 \\3\end{array} \right)| + \frac{1}{2}  |(|\bar a_{n-1}|, |\bar b_{n-1}|)\cdot  \left(\begin{array}{c} 1 \\1\end{array} \right)| = |-|\bar a_{n-1}| + \frac{3}{2} |b_{n-1}|| + |\frac{1}{2} |\bar a_{n-1}| + \frac{1}{2} |\bar b_{n-1}|| \le \frac{3}{2} |a_{n-1}| + 2|b_{n-1}|$.

\medskip

Similarly, $|\bar b_n| = \frac{1}{2} |(|\bar a_{n-1}|, |\bar b_{n-1}|) \cdot \left(\begin{array}{c} 1\\ 1\end{array} \right)| + \frac{1}{2}  |(\bar |a_{n-1}|, |\bar b_{n-1}|)\cdot  \left(\begin{array}{c} 0 \\1\end{array} \right)| = |\frac{1}{2} |\bar a_{n-1}| + \frac{1}{2} |\bar b_{n-1}|| + \frac{1}{2} |\bar b_{n-1}| \le \frac{1}{2} |\bar a_{n-1}| + |\bar b_{n-1}|$.

\medskip

Therefore, we can write $(|\bar a_n|, |\bar b_n|) \le M' \cdot (|\bar a_{n-1}|, |\bar b_{n-1}|)$, meaning the inequality is satisfied by each  coordinate of a vector.
The largest eigenvalue of the matrix $M'=\frac{1}{2}|A|+\frac{1}{2}|B|= \left(\begin{array}{cc} \frac{3}{2} & \frac{1}{2} \\ 2 & 1 \end{array} \right)$ is $1+ \sqrt{2}$, so the average growth rate in this case is $\le 1+ \sqrt{2}$.



\subsection{Proof of Proposition \ref{thm1}}\label{cut}
\medskip

\begin{proof}

\noindent We start by proving part (a).
\medskip

\noindent {\bf (a)} As mentioned above, instead of producing a matrix $M$ from a system of recurrence relations, one can just compute the matrix $M$ of such a system directly, as $M=pA + (1-p)B$. If all eigenvalues of the matrix $M$ are distinct real numbers, then the largest eigenvalue (the Perron–Frobenius eigenvalue) of this matrix will give the average growth rate we are looking for because it will determine the growth rate of the solutions of the corresponding system of recurrence relations.
This establishes part (a) of Proposition \ref{thm1}.
\medskip

\noindent {\bf (b)} What changes here compared to part (a) is that equality like\\
$|\bar a_n| = \frac{1}{2} (|m_1 \cdot |\bar a_{n-1}| + m_2 \cdot |\bar b_{n-1}|) + \frac{1}{2} (|m_3 \cdot |\bar a_{n-1}| + m_4 \cdot |\bar b_{n-1}|)$ (where $m_i$ are coefficients that could be negative) implies $|\bar a_n| \le (|m_1|+|m_3|) |\bar a_{n-1}| + (|m_2| + |m_4|) |\bar b_{n-1}|$ (i.e., inequality instead of equality).

The result follows (cf. the example in Section \ref{growthav1-1}).
\end{proof}

\medskip

\subsection{Average growth for products of $A(k)$ and $B(k)$}\label{growthav22}\hfill \break

For the pair $(A(k), B(k))$ (with positive $k$) with $p=\frac{1}{2}$, we have $M = \frac{1}{2}(A(k)+B(k))=\left(
 \begin{array}{cc} 1 & \frac{k}{2} \\ \frac{k}{2} & 1 \end{array} \right).$

The largest eigenvalue of this matrix is $1+\frac{k}{2}$, so the average growth rate of entries in a random product of matrices $A(k)$ and $B(k)$ is $\mu=1+\frac{k}{2}$.
\medskip

If $k$ is negative, then our Proposition \ref{thm1} only gives that the average growth rate is
$\le 1+\frac{|k|}{2}$. However, in this special case, since $A(-k)=A(k)^{-1}, B(-k)=B(k)^{-1},$
and the determinant of $A(k)$ and $B(k)$ is 1, the average growth rate of entries in a random product of matrices $A(k)$ and $B(k)$ must be equal to that for a random product of $A(-k)$ and $B(-k)$, and is therefore exactly $1+\frac{|k|}{2}$.

\subsection{Average growth for products of matrices from Pollicott's paper}\label{Pollicottmatrices}\hfill \break

\noindent To conclude this section, we treat a pair of matrices from Pollicott's paper \cite{Pollicott}.

In \cite{Pollicott}, Pollicott considered the following two matrices:
$A = \left(
 \begin{array}{cc} 2 & 1 \\ 1 & 1 \end{array} \right) , \hskip .2cm B = \left(
 \begin{array}{cc} 3 & 1 \\ 2 & 1 \end{array} \right)$.

Using the same method as in the previous subsections, we have the following system of recurrence relations for the expectations:
\smallskip

$\bar a_n = \frac{1}{2}(2\bar a_{n-1} + \bar b_{n-1}) + \frac{1}{2}(3\bar a_{n-1} + 2\bar b_{n-1}) = \frac{5}{2}\bar a_{n-1} + \frac{3}{2}\bar b_{n-1}$.

$\bar b_n = \frac{1}{2}(\bar a_{n-1} + \bar b_{n-1}) + \frac{1}{2}(\bar a_{n-1} + \bar b_{n-1}) = \bar a_{n-1} + \bar b_{n-1}$.

Therefore, the corresponding matrix $M$ is $\left( \begin{array}{cc} \frac{5}{2} & \frac{3}{2} \\ 1 & 1 \end{array} \right)$, and the largest eigenvalue of the matrix $M$ is
$\frac{7+\sqrt{33}}{4}$.

Thus, in this case the average growth rate is $\mu=\frac{7+\sqrt{33}}{4} \approx 3.186$.

\section{Upper bounds on the Lyapunov exponent}\label{upper}

Now we are going to look at the growth of the entries in a random product of $n$ matrices $A$ and $B$ of size $d \times d$, where each factor is $A$ or $B$ with probability $p$ and $1-p$, respectively. If the largest entry in such a product is of magnitude $\Theta(s^n)$, then we will call $s$ the {\it generic growth rate} of entries in a random matrix product.

Pollicott \cite{Pollicott} considered a random product of matrices where each matrix comes from a finite collection $\{M_1, \ldots, M_k\}$ of matrices and studied the growth of the {\it norm} of such a product when the number $n$ of factors in a product goes to infinity. To relate what we call the generic growth rate $s$ to what is called the (maximal) Lyapunov exponent $\lambda$ in the theory of stochastic processes, we mention that by a classical  result of Kesten and Furstenberg \cite{Kesten},

\begin{equation}\label{lambda}
\lambda = \lim_{n \to \infty} \frac{1}{n} \log ||M_1 M_2 \cdots M_n|| = \log (\lim_{n \to \infty} \sqrt[n]{||M_1 M_2 \cdots M_n||}),
\end{equation}

\noindent where $\log$ denotes the natural logarithm and $||M||$ denotes the norm of a matrix $M$. Any sub-multiplicative norm can be used here.


Pollicott \cite{Pollicott} offered an algorithm that allows to determine the growth rate of the largest entry in a random product of $n$ {\it non-singular  strictly positive} matrices $A$ and $B$ with any desired precision.

In this section, we offer a simple and efficient method of estimating the Lyapunov exponent based on Theorem \ref{thm2} from the Introduction. First we give a proof of Theorem \ref{thm2}.

\subsection{Proof of Theorem \ref{thm2}}\label{proof2}

\begin{proof}
Denote by $M_n=X_1X_2\cdots X_n$ a random product of $n$ matrices of size $d \times d$, where each $X_i$ is either $A$ or $B$. By $F(n)=||X_1X_2\cdots X_n||$ we denote the {\it max norm} of the matrix $M_n$, i.e., the maximum of the absolute values of the entries of $M_n$. To make this norm sub-multiplicative, we scale it by multiplying it by $d$.

By \cite[Theorem 1, inequality (2.3)]{Kesten}, we have:

\begin{equation}\label{eq1}
\lim_{n \to \infty} \frac{1}{n} \log(F(n)) \le \lim_{n \to \infty} \frac{1}{n} E[\log(F(n))]
\end{equation}

\noindent with probability 1. At the same time, by Jensen's lemma we have

\begin{equation}\label{eq2}
\lim_{n \to \infty} \frac{1}{n} E[\log(F(n))] \le \lim_{n \to \infty} \frac{1}{n} \log(E[(F(n))]) =  \lim_{n \to \infty} \log(\sqrt[n]{E[(F(n)]}) \le  \log \mu.
\end{equation}

\noindent The last inequality follows from our Proposition \ref{thm1}.

\end{proof}

\section{Examples}\label{Examples}

In all examples in this section, we handle products of matrices where each factor is $A$ or $B$ with probability $\frac{1}{2}$.
\medskip

For the pair of matrices $A(k) = \left(
 \begin{array}{cc} 1 & k \\ 0 & 1 \end{array} \right) , \hskip .2cm B(m) = \left(
 \begin{array}{cc} 1 & 0 \\ m & 1 \end{array} \right)$ with positive $k$ and $m$, Corollary 1 from \cite{Sturman} gives the following upper bound for the Lyapunov exponent: $\lambda \le \frac{1}{4}[c + \log(\sqrt{km} + 1/\sqrt{km}) + \frac{1}{2} \log(1+km)]$, where $c$ is a constant approximately equal to 1.0157.

Our method gives $M = \frac{1}{2}\left(
 \begin{array}{cc} 1 & k \\ 0 & 1 \end{array} \right) + \frac{1}{2}\left(
 \begin{array}{cc} 1 & 0 \\ m & 1 \end{array} \right) =\left(\begin{array}{cc} 1 & \frac{k}{2} \\ \frac{m}{2} & 1 \end{array} \right)$, and the largest eigenvalue of the matrix $M$ is $\mu=1+\frac{1}{2} \sqrt{km}$, so  $\log \mu= \log(1+\frac{1}{2} \sqrt{km})$. Note that our method works if either both $k$ and $m$ are positive or both are negative.
\medskip

Thus, if $km$ is large, the method of \cite{Sturman} gives a tighter upper bound for the Lyapunov exponent $\lambda$, whereas our method gives a tighter upper bound for small $km$.
We give a couple of specific instances below.
\medskip

\noindent {\bf (1)} In the case of matrices $A(1)$ and $B(1)$  (see Section \ref{growthav11}), or more generally, $A(k)$ and $B(\frac{1}{k})$, we have $\mu=1.5$, so our upper bound for $\lambda$ in this case is $\log 1.5 \approx 0.405$. To compare, the upper bound provided by Corollary 1 from \cite{Sturman} is 0.514, so our upper bound is significantly tighter in this case.
\medskip

\noindent {\bf (2)} For the matrices $A(2)$ and $B(2)$, Corollary 1 from \cite{Sturman} gives $\lambda \le 0.684$.  We know from Section \ref{growthav22} that $\mu=2$ in this case. Therefore, $\lambda \le \log 2 \approx 0.693$. Thus, for $A(2)$ and $B(2)$ the upper bound provided by Corollary 1 from \cite{Sturman} is better than ours.
\medskip

\noindent {\bf (3)} More generally, for the matrices $A(k)$ and $B(k)$, the upper bound from \cite[Corollary 1]{Sturman} gives $\lambda \le \frac{1}{4}[c + \log(k+\frac{1}{k}) + \frac{1}{2} \log(1+k^2)]$, which is asymptotically equal to $\frac{1}{4}c + \frac{1}{2} \log k$. At the same time, our method (see Section \ref{cut}) gives  $\lambda \le \log(1+\frac{k}{2})$. Thus, our method gives a tighter upper bound on $\lambda$ for small (positive) $k$, whereas the method of \cite{Sturman} gives a tighter upper bound for larger $k$.
\medskip

\noindent {\bf (4)} For Pollicott's matrices $A= \left(
 \begin{array}{cc} 2 & 1 \\ 1 & 1 \end{array} \right) , \hskip .2cm B = \left(
 \begin{array}{cc} 3 & 1 \\ 2 & 1 \end{array} \right)$, we have $\mu=\frac{7+\sqrt{33}}{4} \approx 3.186$, so our upper bound for $\lambda$ in this case is $\log 3.186 \approx 1.159$, whereas Pollicott \cite{Pollicott} gives the following approximation for the actual value of $\lambda$: 1.1433...
\medskip

\noindent {\bf (5)} For the matrices $A= \left(
 \begin{array}{cc} 3 & 1 \\ 1 & 3 \end{array} \right) , \hskip .2cm B = \left(
 \begin{array}{cc} 5 & 2 \\ 2 & 5 \end{array} \right)$ from the paper \cite{Jurga}, we have $\mu=5.5$, so our upper bound for $\lambda$ in this case is $\log 5.5 \approx 1.7$. The upper bound in \cite{Jurga} is 1.66.., so it is tighter.
\medskip

\noindent {\bf (6)} Unlike analytical methods of \cite{Jurga}, \cite{Pollicott2}, and \cite{Sturman}, our method also works for some matrices that have negative entries, as long as the ``expectation matrix" $M$ has only positive entries, or more generally, satisfies the conditions of our Proposition \ref{thm1}. For example, let $A= \left(
 \begin{array}{cc} 1 & -1 \\ 0 & 1 \end{array} \right), \hskip .2cm B = \left(
 \begin{array}{cc} 1 & 2 \\ 2 & 1 \end{array} \right)$. Then $M = \frac{1}{2}A + \frac{1}{2}B = \left(\begin{array}{cc} 1 & \frac{1}{2} \\ 1 & 1 \end{array} \right)$, so the largest eigenvalue of the matrix $M$ is $\mu=1+ \sqrt{\frac{1}{2}} \approx 1.707$,   and $\log \mu \approx 0.535$, which is an upper bound on the Lyapunov exponent in this case.
\medskip

\subsection{Series of matrices from Pollicott's paper \cite{Pollicott2}}\hfill \break

\noindent In \cite{Pollicott2}, Pollicott treated the following series of pairs of matrices:

$A_1(t)= \left(
 \begin{array}{cc} 1+t & 1 \\ t & 1 \end{array} \right), \hskip .2cm A_2(t) = \left(
 \begin{array}{cc} 1 & t \\ 1 & 1+t \end{array} \right)$, where $t>0$.

The ``expectation matrix" for the Bernoulli distribution $(\frac{1}{2}, \frac{1}{2})$ in this case is $M(t) = \left(\begin{array}{cc} 1+\frac{t}{2} & \frac{1}{2} + \frac{t}{2} \\ \frac{1}{2} + \frac{t}{2} & 1+\frac{t}{2} \end{array} \right)$. The largest eigenvalue of the matrix $M(t)$ is $\mu=t+1.5$, so our upper bound for the Lyapunov exponent in this case is $\log(t+1.5)$.
\medskip

For particular values of $t$ from \cite{Pollicott2}, we have the following comparisons of  upper bounds:
\medskip

$\bullet$ $t=2:$  $\log \mu \approx 1.2528$, whereas the upper bound from \cite{Pollicott2} is $\approx 1.2509$.
\medskip

$\bullet$ $t=1:$  $\log \mu \approx 0.916$, whereas the upper bound from \cite{Pollicott2} is $\approx 0.915$.
\medskip

$\bullet$ $t=0.5:$  $\log \mu \approx 0.6931$, whereas the upper bound from \cite{Pollicott2} is $\approx 0.6936$.
\medskip

$\bullet$ $t=0.4:$  $\log \mu \approx 0.6418$, whereas the upper bound from \cite{Pollicott2} is $\approx 0.6468$.

\medskip

$\bullet$ $t=0.3:$  $\log \mu \approx 0.5878$, whereas the upper bound from \cite{Pollicott2} is $\approx 0.5872$.
\medskip

$\bullet$ $t=0.2:$  $\log \mu \approx 0.5306$, whereas the upper bound from \cite{Pollicott2} is $\approx 0.529$.
\medskip

$\bullet$ $t=0.1:$  $\log \mu \approx 0.47$, whereas the upper bound from \cite{Pollicott2} is $\approx 0.4$.
\medskip

Thus, we have an interesting behavior here: while our upper bounds are very close to those of \cite{Pollicott2} (with the exception of the case $t=0.1$), for $t=0.4$ and $t=0.5$ our bounds are tighter, whereas for other values of $t$ above Pollicott's bounds are tighter.

\section{Conclusions}\label{Conclusions}

To summarize, our combinatorial method of producing an upper bound for the Lyapunov exponent of a random product of matrices has both advantages and disadvantages compared to analytical methods of \cite{Jurga}, \cite{Pollicott}, \cite{Pollicott2}, \cite{Sturman}, and others.

One obvious advantage of our method is unparalleled simplicity and efficiency. Another advantage is applicability of our method to a larger class of matrices. Most of the time, methods such as those in \cite{Jurga}, \cite{Pollicott}, \cite{Pollicott2}, \cite{Sturman} require participating matrices to be non-singular and have only positive entries.
Our method is insensitive to that.

To be fair, several authors, e.g. \cite{Protasov}, \cite{Sturman}, managed to expand the applicability of their methods to include some matrices with negative entries (\cite[Section 2.3]{Sturman}) and some singular matrices (\cite{Protasov}), although the class of matrices covered by our method still appears to be larger. We regret to be unable to offer any insights into why this is the case; the methods are just too different, to the point of having nothing in common.

On the other hand, an obvious disadvantage of our method is that it cannot provide any error estimates of the obtained bounds on the Lyapunov exponent.

As for tightness of upper bounds, we have seen in the previous sections that
our upper bounds are sometimes better, sometimes not, compared to those of the above cited papers. 

It is hardly possible to tell when or why our method produces tighter upper bounds because, again, our combinatorial method is way too different from the existing analytical methods. Given any particular collection of matrices, the easiest way to compare the corresponding upper  bounds for the Lyapunov exponent is just to compute them using different methods and compare the results.

In any case, it is always good to have different approaches  to the same problem.

\end{document}